\documentclass[12pt]{article}
\usepackage{amsfonts}
\usepackage{theorem}
\newtheorem{theorem}{Theorem}
\newtheorem{corollary}[theorem]{Corollary}
\usepackage[unicode,colorlinks,plainpages=false]{hyperref}
\newcommand{\openbox}{$\begin{array}{c}
\hspace*{-0.55em}\sqcap \hspace*{-0.60em}\\[-0.4em] \hline
\multicolumn{1}{c}{\hspace*{-0.60em}}\\[-0.8em]
\end{array}$}

\begin{document}

\centerline{\bf On Monoid Congruences of Commutative}
\centerline{\bf Semigroups}

\bigskip

\bigskip

\centerline{Attila Nagy}

\bigskip

\bigskip

\begin{abstract}
Let $S$ be a semigroup and $A$ a subset of $S$. By the separator $SepA$ of $A$ we mean the set of all elements $x\in S$ which satisfy $xA\subseteq A$, $Ax\subseteq A$, $x(S\setminus A)\subseteq (S\setminus A)$, $(S\setminus A)x\subseteq (S\setminus A)$. In this paper we characterize the monoid congruences of commutative semigroups by the help of the notion of the separator. We show that every monoid congruence of a commutative semigroup $S$ can be constructed by the help of subsets $A$ of $S$ for which $SepA\neq \emptyset$.
\end{abstract}

\bigskip

Let $S$ be a semigroup and $A$ a subset of $S$. By the idealizer of $A$ we mean the set of all elements $x$ of $S$ which satisfy $xA\subseteq A$ and $Ax\subseteq A$. The idealizer of $A$ will be denoted by $IdA$. As in \cite{2}, $IdA\cap Id(S\setminus A)$ is called the separator of $A$ and will be denoted by $SepA$.

In this paper we characterize the monoid congruences of commutative semigroups by the help of the separator. We show that a commutative semigroup $S$ has a non universal monoid congruence if and only if $SepA\neq \emptyset$ for some subset $A$ of $S$ with $\emptyset \subset A\subset S$. Moreover, every monoid congruence on a commutative semigroup $S$ can be constructed by the help of subsets $A$ of $S$ for which $SepA\neq \emptyset$.

\medskip

\underline{Notations}. Let $S$ be a semigroup and $H$ a subset of $S$. Following \cite{1}, let
\[H\dots a=\{(x,y)\in S\times S:\ xay\in H\}, \quad a\in S\]
and
\[P_H=\{(a,b)\in S\times S:\ H\dots a=H\dots b\}.\]
If $\{H_i, i\in I\}$ is a family of subsets $H_i$ of $S$ such that $H=\cap _{i\in I}SepH_i$, then the family $\{H_i, i\in I\}$ will be denoted by $(H;H_i,I)$.
For a family $(H;H_i,I)\neq \emptyset$, we define a relation $P(H;H_i,I)$ on $S$ as follows:
\[P(H;H_i,I)=\{(a,b)\in S\times S:\ H_i\dots a=H_i\dots b \quad \hbox{for all}\quad i\in I\}.\]
For notations and notions not defined here, we refer to \cite{1} and \cite{2}.

\newpage

\begin{theorem}\label{th1} Let $S$ be a semigroup and $p$ a congruence on $S$. If $S_k$ ($k\in K$) is a family of congruence classes of $S$ modulo $p$, then the separator of $\cup _{k\in K}S_k$ is either empty or the union of some congruence classes of $S$ modulo $p$.
\end{theorem}

\noindent
{\bf Proof}. Let $S_k$ ($k\in K$) be a family of congruence classes of $S$ modulo $p$, and let $U=\cup _{k\in K}S_k$. We may assume $SepU\neq \emptyset$ and $SepU\neq S$. Then there exist elements $a, b\in S$ such that $a\in SepU$ and $b\notin SepU$. We consider an arbitrary couple $(a, b)$ with this property, and prove that $(a, b)\notin p$. By the assumption, at least one of the following four condition holds for $b$:

\begin{enumerate}
\item[(1.1)] $bU\not\subseteq U$,
\item[(1.2)] $Ub\not\subseteq U$,
\item[(1.3)] $b(S\setminus U)\not\subseteq (S\setminus U)$,
\item[(1.4)] $(S\setminus U)b\not\subseteq (S\setminus U)$.
\end{enumerate}
In case (1,1), there exists an element $c\in U$ such that $bc\notin U$. Thus $abc\notin U$, because $a\in SepU$. Since $SepU$ is a subsemigroup of $S$ and $c\in U$, we have $aac\in U$. As $U$ is the union of congruence classes of $S$ modulo $p$, our result implies that $a$ and $b$ do not belong to the same congruence class of $S$ modulo $p$. The same conclusion holds in cases (1.2), (1.3) and (1.4), too. From this it follows that $SepU$ is the union of congruence classes of $S$ modulo $p$.\hfill\openbox

\begin{theorem}\label{th2} Let $S$ be a semigroup and $H$ a subsemigroup of $S$. If $(H;H_i,I)$ is a non empty family of subsets of $S$, then $P(H;H_i,I)$ is a congruence on $S$ such that the subsets $H_i$ ($i\in I$) and $H$ are unions of some congruence classes of $S$ modulo $P(H;H_i,I)$.
\end{theorem}

\noindent
{\bf Proof}. It can be easily verified that $P(H;H_i, I)$ is a congruence on $S$. Let $i\in I$ be abitrary. Assume $H_i\neq S$. Let $x, y\in S$ such that $x\in H_i$, $y\notin H_i$. Let $h\in H$. Since $H\subseteq SepH_i$, we have $hxh\in H_i$ and $hyh\notin H_i$. Thus $(x, y)\notin P(H; H_i, I)$ and so $H_i$ is the union of some congruence classes of $S$ modulo $P(H; H_i, I)$.

To show that $H$ is the union of some congruence classes of $S$ modulo $P(H; H_i, I)$ let $h\in H$ and $g\in (S\setminus H)$ be arbitrary elements. Then there is an index $j$ in $I$ such that $g\notin SepH_j$. From this it follows that at least one of the following holds for $g$:

\newpage

\begin{enumerate}
\item[(1.5)] $gH_j\not\subseteq H_j$,
\item[(1.6)] $H_jg\not\subseteq H_j$,
\item[(1.7)] $g(S\setminus H_j)\not\subseteq (S\setminus H_j)$,
\item[(1.8)] $(S\setminus H_j)g\not\subseteq (S\setminus H_j)$.
\end{enumerate}
In case (1.5), there exists an element $b$ in $H_j$ such that $gb\notin H_j$. Then $hgb\notin H_j$. As $hhb\in H_j$, we have $(g, h)\notin P(H;H_i, I)$. The same conclusion holds in cases (1.6), (1.7) and (1.8), too. Consequently, $H$ is the union of some congruence classes of $S$ modulo $P(H;H_i, I)$. Thus the theorem is proved.\hfill\openbox

\begin{theorem}\label{th3} Let $S$ be a commutative semigroup and $H$ a subsemigroup of $S$. Assume that $(H;H_i,I)$ is a non empty family of subsets of $S$. Then $P(H;H_i,I)$ is a monoid congruence on $S$ such that $H$ is the identity element of $S/P(H;H_i,I)$.
Conversely, every monoid congruence on a commutative semigroup can be so constructed.
\end{theorem}

\noindent
{\bf Proof}. Let $S$ be a ommutative semigroup and $H$ a subsemigroup of $S$. Assume that $(H; H_i, I)$ is not empty. Then, by Theorem~\ref{th2}, $H$ is a union of some congruence classes of $S$ modulo $P(H; H_i, I)$. Let $a, b\in H$. We show that $(a, b)\in P(H; H_i, I)$. Let $i\in I$ and $x, y\in S$ be arbitrary. Assume $xay\in H_i$. Then $yxa\in H_i$ and so $yx\in H_i$, because $S$ is commutative and $a\in H\subseteq SepH_i$. Thus $yxb\in H_i$ and so $xby\in H_i$, because $b\in H\subseteq SepH_i$. We can prove similarly that $xay\notin H_i$ implies $xby\notin H_i$. Thus $(a, b)\in P(H; H_i, I)$, indeed. Consequently, $H$ is a congruence class of $S$ modulo $P(H; H_i, I)$.

Next we show that $H$ is the identity element of the factor semigroup $S/P(H; H_i, I)$. Let $S_k$ be an arbitrary congruence class of $S$ modulo \break
$P(H; H_i, I)$. Let $u\in S_k$ be arbitrary. We show that, for any $a\in H$, the product $ua$ belongs to $S_k$. Let $i\in I$ and $x, y\in S$ be arbitrary elements. Since $S$ is commutative and $a\in H\subseteq SepH_i$, the product $xuy$ belongs to $H_i$ if and only if $xuay=xuya$ belongs to $H_i$. Thus $(u, ua)\in P(H; H_i, I)$ and so $ua\in S_k$. Thus $H$ is the identity element of the factor semigroup $S/P(H; H_i, I)$, indeed.

Conversely, let $S$ be a commutative semigroup and $p$ a monoid congruence on $S$. Denote $H$ the identity element of the factor semigroup $S/p$. Let $M=\cap _{k\in K}SepS_k$, where $\{S_k, k\in K\}$ is the set of all congruence classes of $S$ modulo $p$. It is clear that $H\subseteq M$. We show that $H=M$. Assume, in an indirect way, that $H\subset M$. Let $a\in H$ and $b\in M\setminus H$ be arbitrary elements. Then there is an element $k_0\in K$ such that $b\in S_{k_0}$. As $b\in M\subseteq SepS_{k_0}$, we have $SepS_{k_0}\cap S_{k_0}\neq \emptyset$ and so $SepS_{k_0}\subseteq S_{k_0}$ (see Theorem 3 of \cite{2}). From this it follows that  $H\subseteq M\subset SepS_{k_0}\subseteq S_{k_0}$ and so $H=S_{k_0}$, because $H$ and $S_{k_0}$ are congruence classes of $S$ modulo $p$. As $b\in S_{k_0}$, we get $b\in H$ which is a contradiction. Hence $H=M$. Consequently the congruence $P(H; S_k, K)$ is defined.

We show that $P(H; S_k, K)=p$. To show  $P(H; S_k, K)\subseteq p$, let $a, b\in S$ be arbitrary elements with $(a, b)\in P(H; S_k, K)$. Let $m, n\in K$ such that $a\in S_m$, $b\in S_n$. Since $H$ is the identity element of the factor semigroup $S/p$, $hah\in S_m$ and $hbh\in S_n$ for an arbitrary $h\in H$. If $n\neq m$ then $(h, h)\in S_m...a$ and $(h, h)\notin S_m...b$, because $hbh\notin S_m$. In this case $(a, b)\notin P(H; S_k, K)$ which is a contradiction. Thus $n=m$ and so $a, b\in S_m=S_n$. Consequently $(a, b)\in p$. Hence $P(H;S_k, K)\subseteq p$.
As $(a, b)\in p$ implies $(xay, xby)\in p$ for all $x, y\in S$, we get $S_k...a=S_k...b$ for all $k\in K$ which implies that $(a, b)\in P(H; S_k, K)$. Consequently $p\subseteq P(H; S_k, K)$. Therefore $p=P(H; S_k, K)$.\hfill\openbox

\medskip

A subset $U$ of a semigroup $S$ is called an unitary subset of $S$ if, for every $a, b\in S$, the assumption $ab,b\in U$ implies $b\in U$, and also
$ab, a\in U$ implies $b\in U$.

\begin{theorem}\label{th4} Let $S$ be a commutative semigroup and $H$ a subsemigroup of $S$. If $p$ is a monoid congruence on $S$ such that $H$ is the identity of $S/p$, then $P(H;H_i,I)\subseteq p\subseteq P_H$, where $\{H_i, i\in I\}$ denotes the family of all subsets $H_i$ of $S$ satisfying $H\subseteq SepH_i$ ($i\in I$).
\end{theorem}

\noindent
{\bf Proof}. Let $p$ be a monoid congruence on a commutative semigroup $S$, and let $H\subseteq S$ be the identity element of $S/p$. Then $H$ is an unitary subsemigroup of $S$ and so $H=SepH$ (see Theorem 8 of \cite{2}). From this it follows that $H=\cap _{i\in I}SepH_i$, where $\{ H_i, i\in I\}$ is the family of all subsets $H_i$ of $S$ for which $H\subseteq SepH_i$. Thus the congruence $P(H; H_i, I)$ is defined on $S$. Let $\{ S_k, k\in K\}$ be the family of all congruence classes of $S$ modulo $p$. By Theorem~\ref{th3}, $p=P(H; S_k, K)$. As $H\in (H; S_k, K)\subseteq (H; H_i, I)$, we have
$P(H;H_i,I)\subseteq p\subseteq P_H$.\hfill\openbox

\begin{corollary}\label{cd5} A commutative semigroup $S$ has a non universal monoid congruence if and only if it has a subset $A$ with $\emptyset \subset A\subset S$ such that $SepA\neq \emptyset$.
\end{corollary}

\noindent
{\bf Proof}. Let $p$ be a non universal monoid congruence on a commutative semigroup $S$ and $A$ the congruence class of $S$ modulo $p$ which is the identity element of the factor semigroup $S/p$. Then $\emptyset \subset A\subset S$. As $A\subseteq SepA$, we have $SepA\neq \emptyset$.

Conversely, let $A$ be a subset of a commutative semigroup $S$ such that $\emptyset \subset A\subset S$ and $SepA \neq \emptyset$. As $SepA\subseteq A$ or $SepA\subseteq (S\setminus A)$ by Theorem 3 of \cite{2}, we have $SepA\neq S$. By Theorem~\ref{th3} of this paper, $SepA$ is the identity element of the factor semigroup $S/P_A$ and so $P_A$ is a non universal monoid congruence on $S$.\hfill\openbox

\bigskip

\noindent
Attila Nagy

\noindent
Department of Algebra

\noindent
Mathematical Institute

\noindent
Budapest University of Technology and Economics

\noindent
e-mail: nagyat@math.bme.hu

\end{document}